\documentclass[oneside,english]{amsart}
\usepackage[T1]{fontenc}
\usepackage[latin9]{inputenc}
\usepackage{amsthm}
\usepackage{amssymb}
\usepackage{esint}

\numberwithin{equation}{section} 
\numberwithin{figure}{section} 
\theoremstyle{plain}
\theoremstyle{plain}
\newtheorem{thm}{Theorem}
  \theoremstyle{plain}
  \newtheorem{cor}[thm]{Corollary}
  \theoremstyle{remark}
  \newtheorem{rem}[thm]{Remark}
  \theoremstyle{remark}
  \newtheorem*{acknowledgement*}{Acknowledgement}

\usepackage{babel}

\begin{document}

\title{On Generalized Hilbert Matrices}

\author{Ruiming Zhang}
\begin{abstract}
In this note, we present a systematic method to explicitly compute
the determinants and inverses for some generalized Hilbert matrices
associated with orthogonal systems with explicit representations.
We expressed the determinant, the inverse and a lower bound for the
smallest eigenvalue of such matrix in terms of the orthogonal system.
\end{abstract}


\curraddr{School of Mathematical Sciences\\
Guangxi Normal University\\
Guilin City, Guangxi 541004\\
P. R. China.}


\email{ruimingzhang@yahoo.com}

\maketitle

\section{Introduction}

The Hilbert matrices are the moment matrices associated with Legendre
polynomials. The generalized Hilbert matrices are the generalized
moment matrices associated with certain orthogonal systems. It may
be interesting to find the exact formulas for their determinants and
inverses. In this note we provide a systematic method from the theory
of orthogonal polynomials to find these formulas through the explicit
representations of their related orthogonal systems. We will demonstrate
that, once the orthogonal system is explicitly known, then the inverse
and determinant of each generalized moment matrix are known explicitly.
Furthermore, we also know a lower bound for the smallest eigenvalue.
Since an orthogonal system may be derived in many different ways,
this method could be very handy. We also present four examples to
show how to apply this method in various situations.

\section{Main Results\label{sec:Main-Results}}

Let $E$ be a complex inner product space with a sequence of linearly
independent vectors $\left\{ u_{n}\right\} _{n=0}^{\infty}$. For
each nonnegative integer $n$, the following matrix is positive definite,\begin{align*}
G_{n} & =\left(m_{j,k}\right)_{j,k=0}^{n},\end{align*}
where\begin{align*}
m_{j,k} & =(u_{j},u_{k}).\end{align*}
There is a unique orthonormal system $\left\{ p_{k}\right\} _{k=0}^{\infty}$
with $p_{n}$ having positive leading coefficient in $u_{n}$, which
could be found through the Gram-Schmidt orthogonalization process.
Each $p_{n}$ is given explicitly by \cite{Ismail} \begin{align*}
p_{n} & =\frac{1}{\sqrt{\det G_{n}\det G_{n-1}}}\det\left(\begin{array}{ccccc}
m_{0,0} & m_{0,1} & m_{0,2} & \dots & m_{0,n}\\
m_{1,0} & m_{1,1} & m_{1,2} & \dots & m_{1,n}\\
\vdots & \vdots & \vdots & \ddots & \vdots\\
m_{n-1,0} & m_{n-1,1} & m_{n-1,2} & \dots & m_{n-1,n}\\
u_{0} & u_{1} & u_{2} & \dots & u_{n}\end{array}\right).\end{align*}
Clearly,\begin{align*}
p_{n} & =\sum_{j=0}^{n}a_{n,j}u_{j},\quad a_{n,n}=\sqrt{\frac{\det G_{n-1}}{\det G_{n}}}.\end{align*}
Thus,\begin{align*}
\det G_{n} & =\prod_{j=0}^{n}a_{j,j}^{-2}.\end{align*}
Let $\left\{ v_{n}\right\} _{n=0}^{\infty}$ be another sequence of
vectors in $E$ related to $\left\{ u_{n}\right\} _{n=0}^{\infty}$
in the following way,\begin{align*}
v_{n} & =\sum_{j=0}^{n}c_{n,j}u_{j},\quad c_{n,n}\neq0,\quad n\ge0,\end{align*}
then, they are also linearly independent. Let us define 

\begin{align*}
H_{n} & =\left(I_{j,k}\right)_{j,k=0}^{n},\quad I_{j,k}=(u_{j},v_{k}),\end{align*}
 then\begin{align*}
H_{n} & =G_{n}C_{n}^{*},\end{align*}
where $C_{n}^{*}$ is the Hermitian conjugate of $C_{n}$ and \begin{align*}
C_{n} & =\left(c_{j,k}\right)_{j,k=0}^{n},\end{align*}
with\[
c_{j,k}=0,\quad k>j.\]
Thus,\begin{align*}
\det H_{n} & =\det G_{n}\prod_{j=0}^{n}\overline{c_{j,j}}.\end{align*}
Let\begin{align*}
p_{n} & =\sum_{k=0}^{n}b_{n,k}v_{k},\end{align*}
and\begin{align*}
A_{n} & =\left(a_{j,k}\right)_{j,k=0}^{n},\quad B_{n}=\left(b_{j,k}\right)_{j,k=0}^{n},\end{align*}
where we follow the same convention as above,\[
a_{j,k}=b_{j,k}=0,\quad k>j.\]
Evidently, both $A_{n}$ and $B_{n}$ are invertible. Let \begin{align*}
A_{n}^{-1} & =\left(s_{j,k}\right)_{j,k=0}^{n},\quad B_{n}^{-1}=\left(t_{j,k}\right)_{j,k=0}^{n},\end{align*}
then,\begin{align*}
u_{j} & =\sum_{\ell=0}^{n}s_{j,\ell}p_{\ell},\quad v_{k}=\sum_{m=0}^{n}t_{k,m}p_{m},\end{align*}
and \begin{align*}
I_{j,k} & =(u_{j},v_{k})=\sum_{m=0}^{n}s_{j,m}\overline{t_{k,m}}.\end{align*}
Thus,\begin{align*}
H_{n} & =A_{n}^{-1}(B_{n}^{-1})^{*}=A_{n}^{-1}\left(B_{n}^{*}\right)^{-1},\end{align*}
which gives \begin{align*}
H_{n}^{-1} & =B_{n}^{*}A_{n}.\end{align*}
We summarize our above discussion as the following theorem:
\begin{thm}
\label{thm:main}For each nonnegative integer $n$, assume that $G_{n}$,
$H_{n}$, $A_{n}$, $B_{n}$ and $C_{n}$ as defined above, then,\begin{align*}
H_{n} & =G_{n}C_{n}^{*},\quad H_{n}^{-1}=B_{n}^{*}A_{n}.\end{align*}

\begin{align*}
\det G_{n} & =\prod_{j=0}^{n}a_{j,j}^{-2},\ \det C_{n}=\prod_{j=0}^{n}c_{j,j},\ \det H_{n}=\prod_{j=0}^{n}\frac{\overline{c_{j,j}}}{a_{j,j}^{2}}\end{align*}
where $A_{n}^{*}$ is the Hermitian conjugate of $A_{n}$. 
\end{thm}
Notice that\begin{align*}
H_{n}^{-1}: & =\left(\gamma_{j,k}\right)_{j,k=0}^{n}\end{align*}
with\begin{align*}
\gamma_{j,k} & =\sum_{\ell=\max(j,k)}^{n}\overline{b_{\ell,j}}a_{\ell,k}.\end{align*}
From the expressions of $\det G_{n}$ and $\gamma_{j,k}$, we observe
that in the case $u_{k}=v_{k}$, the requirement that $a_{k,k}>0$
could be disregarded in the actual computations, for, if we replace
the pair $\left\{ u_{k},v_{k}\right\} $ by $\left\{ \epsilon_{k}u_{k},\epsilon_{k}v_{k}\right\} $
with $|\epsilon_{k}|=1$ we won't change $H_{n}$. In some situations,
the entries for $G_{n}$ are too complicated to be useful. This is
the reason we won't compute $G_{n}$ for the Askey-Wilson polynomials,
we compute $H_{n}$ instead. But $G_{n}$ could be recovered easily
via\begin{align*}
G_{n} & =H_{n}\left(C_{n}^{*}\right)^{-1},\quad\det G_{n}=\frac{\det H_{n}}{\prod_{j=0}^{n}\overline{c_{j,j}}}.\end{align*}

Let $\lambda_{s}$ be the smallest eigenvalue of $G_{n}$ , then $\frac{1}{\lambda_{s}}$
is the largest eigenvalue of $G_{n}^{-1}$. Since $G_{n}^{-1}$ is
positive definite, we have\[
\frac{1}{\lambda_{s}}=\left\Vert G_{n}^{-1}\right\Vert _{2}\le\sqrt{\sum_{i,j=0}^{n}\left|\rho_{i,j}\right|^{2}}\]
with\[
G_{n}^{-1}=\left(\rho_{i,j}\right)_{i,j=0}^{n},\]
and\begin{align*}
\rho_{j,k} & =\sum_{\ell=\max(j,k)}^{n}\overline{a_{\ell,j}}a_{\ell,k}.\end{align*}
Hence,\[
\left|\rho_{j,k}\right|^{2}\le\sum_{\ell=j}^{n}\left|a_{\ell,j}\right|^{2}\sum_{\ell=k}^{n}\left|a_{\ell,k}\right|^{2},\]
and\[
\frac{1}{\lambda_{s}}\le\sum_{j=0}^{n}\sum_{\ell=j}^{n}\left|a_{\ell,j}\right|^{2}.\]

Another lower bound could be found by considering the $||\cdot||_{\infty}$.
Let $\left(x_{0},\dots,x_{n}\right)^{T}$ be an eigenvector of $G_{n}^{-1}$
corresponding to $\lambda_{s}^{-1}$, then we have\[
\frac{1}{\lambda_{s}}\left(\begin{array}{c}
x_{0}\\
\vdots\\
x_{n}\end{array}\right)=\left(\begin{array}{ccc}
\rho_{0,0} & \dots & \rho_{0,n}\\
\vdots & \vdots & \vdots\\
\rho_{n,0} & \dots & \rho_{n,n}\end{array}\right)\left(\begin{array}{c}
x_{0}\\
\vdots\\
x_{n}\end{array}\right)=\left(\begin{array}{c}
{\displaystyle \sum_{k=0}^{n}}\rho_{0,k}x_{k}\\
\vdots\\
{\displaystyle \sum_{k=0}^{n}}\rho_{n,k}x_{k}\end{array}\right),\]
 then,\[
\frac{1}{\lambda_{s}}\left\Vert \left(\begin{array}{c}
x_{0}\\
\vdots\\
x_{n}\end{array}\right)\right\Vert _{\infty}=\left\Vert \left(\begin{array}{c}
{\displaystyle \sum_{k=0}^{n}}\rho_{0,k}x_{k}\\
\vdots\\
{\displaystyle \sum_{k=0}^{n}}\rho_{n,k}x_{k}\end{array}\right)\right\Vert _{\infty}\le\max_{0\le j\le n}\left\{ {\displaystyle \sum_{k=0}^{n}}\left|\rho_{j,k}\right|\right\} \left\Vert \left(\begin{array}{c}
x_{0}\\
\vdots\\
x_{n}\end{array}\right)\right\Vert _{\infty}.\]
 Observe that\[
{\displaystyle \sum_{k=0}^{n}}\left|\rho_{j,k}\right|\le{\displaystyle \sum_{k=0}^{n}}\sum_{\ell=0}^{n}\left|a_{\ell,j}\right|\left|a_{\ell,k}\right|=\sum_{\ell=0}^{n}\left|a_{\ell,j}\right|{\displaystyle \sum_{k=0}^{n}}\left|a_{\ell,k}\right|,\]
 and\begin{eqnarray*}
\max_{0\le j\le n}\left\{ {\displaystyle \sum_{k=0}^{n}}\left|\rho_{j,k}\right|\right\}  & \le & \max_{0\le j\le n}\left\{ \sum_{\ell=0}^{n}\left|a_{\ell,j}\right|{\displaystyle \sum_{k=0}^{n}}\left|a_{\ell,k}\right|\right\} \\
 & \le & \sum_{\ell=0}^{n}\max_{0\le j\le n}\left\{ \left|a_{\ell,j}\right|\right\} {\displaystyle \sum_{k=0}^{n}}\left|a_{\ell,k}\right|\\
 & \le & \sum_{\ell=0}^{n}\left\{ \sum_{j=0}^{n}\left|a_{\ell,j}\right|{\displaystyle \cdot\sum_{k=0}^{n}}\left|a_{\ell,k}\right|\right\} \\
 & = & \sum_{\ell=0}^{n}\left\{ \sum_{j=0}^{\ell}\left|a_{\ell,j}\right|{\displaystyle \cdot\sum_{k=0}^{\ell}}\left|a_{\ell,k}\right|\right\} \\
 & = & \sum_{\ell=0}^{n}\left\{ \sum_{j=0}^{\ell}\left|a_{\ell,j}\right|\right\} ^{2}.\end{eqnarray*}

\begin{thm}
\label{thm:bounds}For each nonnegative integer $n$, assume that
$G_{n}$ and $A_{n}$, as defined above, then,\[
\lambda_{s}\ge\max\left\{ \frac{1}{\sum_{\ell=0}^{n}\sum_{j=0}^{\ell}\left|a_{\ell,j}\right|^{2}},\frac{1}{\sum_{\ell=0}^{n}\left\{ \sum_{j=0}^{\ell}\left|a_{\ell,j}\right|\right\} ^{2}}\right\} ,\]
where $\lambda_{s}$ is the smallest eigenvalue of $G_{n}$. 
\end{thm}
The second lower bound is particular interesting when the generalized
orthogonal system are certain orthogonal polynomials, since it could
be expressed in terms of the orthonormal polynomials.
\begin{cor}
Let\[
p_{n}(x)=\sum_{k=0}^{n}a_{n,k}x^{k},\quad n=0,1,\dots\]
be the orthonormal polynomials and $m_{n}$ be the $n$-th power moment
with respect to a probability measure $d\mu(x)$. If there is a complex
number $z_{0}$ with $|z_{0}|=1$ such that all of \[
a_{n,k}z_{0}^{k},\quad k=0,1,\dots\]
are of the same sign, then the smallest eigenvalue $\lambda_{s}$
of the matrix\[
G_{n}=\left(m_{i+j}\right)_{i,j=0}^{n}\]
has a lower bond\[
\lambda_{s}\ge\frac{1}{\sum_{\ell=0}^{n}|p_{\ell}(z_{0})|^{2}}.\]

\end{cor}
The first observation for the last corollary is that it seems such
$z_{0}$ always exits. For the symmetric orthogonal polynomials, from
the three term recurrence we see that $z_{0}=i$. But we don't know
how to prove this for general cases. The second observation is that
if all the polynomials are real, then the special case of Christoffel-Darboux
formula gives\[
\sum_{m=0}^{n}|p_{m}(z_{0})|^{2}=\frac{a_{n,n}\left\{ p'_{n+1}(z_{0})p_{n}(z_{0})-p'_{n}(z_{0})p_{n+1}(z_{0})\right\} }{a_{n+1,n+1}},\]
which gives us \[
\lambda_{s}\ge\frac{a_{n+1,n+1}}{a_{n,n}\left\{ p'_{n+1}(z_{0})p_{n}(z_{0})-p'_{n}(z_{0})p_{n+1}(z_{0})\right\} },\]
and this may be useful to find the asymptotic behaviour of the lower
bound.

\subsection{Matrices associated M\" untz systems}
\begin{thm}
\label{thm:muntz1}For $n\in\mathbb{N}$ and $\left\{ \alpha_{0},\alpha_{1},\dots\right\} \subset\mathbb{C}$,
the matrix\begin{align}
 & \left(\frac{1}{\alpha_{j}+\overline{\alpha_{k}}+1}\right)_{j,k=0}^{n}\label{eq:muntz matrix}\end{align}
has determinant\begin{align*}
\det\left(\frac{1}{\alpha_{j}+\overline{\alpha_{k}}+1}\right)_{j,k=0}^{n} & ={\displaystyle \prod_{k=0}^{n}}\frac{{\displaystyle \prod_{j=0}^{k-1}}|\alpha_{k}-\alpha_{j}|^{2}}{(1+2\Re(\alpha_{k})){\displaystyle \prod_{j=0}^{k-1}}|\alpha_{k}+\overline{\alpha_{j}}+1|^{2}}.\end{align*}
Under the condition \begin{align*}
\alpha_{j}\neq\alpha_{k},\,\Re(\alpha_{j})+\Re(\alpha_{k})\neq-1,\, & j,k=0,1,\dots n,\end{align*}
the matrix \eqref{eq:muntz matrix} is invertible, and its inverse
matrix $(\gamma_{j,k})_{j,k=0}^{n}$ has element\begin{align*}
\gamma_{j,k} & =\sum_{m=\max(j,k)}^{n}\frac{(1+2\Re(\alpha_{m})){\displaystyle {\displaystyle \prod_{r=0}^{m-1}}(\overline{\alpha_{j}}+\alpha_{r}+1)(\alpha_{k}+\overline{\alpha_{r}}+1)}}{{\displaystyle \prod_{\begin{array}{c}
p=0\\
p\neq j\end{array}}^{m-1}}(\overline{\alpha_{j}}-\overline{\alpha_{p}}){\displaystyle \prod_{\begin{array}{c}
q=0\\
q\neq k\end{array}}^{m-1}}(\alpha_{k}-\alpha_{q})}.\end{align*}
 When the matrix \eqref{eq:muntz matrix} is positive definite, its
smallest eigenvalue has a lower bound\textup{\begin{align*}
\lambda_{s} & \ge\left\{ \sum_{\ell=0}^{n}\left\{ \sum_{j=0}^{\ell}\frac{\sqrt{1+2\Re(\alpha_{\ell})}{\displaystyle \prod_{k=0}^{\ell-1}|\alpha_{j}+\overline{\alpha_{k}}+1|}}{{\displaystyle \prod_{\begin{array}{c}
k=0\\
k\neq j\end{array}}^{\ell-1}}|\alpha_{j}-\alpha_{k}|}\right\} ^{2}\right\} ^{-1}.\end{align*}
}
\end{thm}
More generalized matrices associated with a generalized M\" untz
system:
\begin{thm}
\label{thm:muntz2}Given $a,b,c\in\mathbb{R}$ and distinct complex
numbers $\left\{ \alpha_{k}\right\} _{k=0}^{\infty}$ , for each positive
integer $n$, the matrix\begin{align}
 & \left(\frac{1}{c\alpha_{j}\overline{\alpha_{k}}-a(\alpha_{j}+\overline{\alpha_{k}})-b}\right)_{j,k=0}^{n}\label{eq:generalized muntz matrix}\end{align}
 has\begin{align*}
 & \det\left(\frac{1}{c\alpha_{j}\overline{\alpha_{k}}-a(\alpha_{j}+\overline{\alpha_{k}})-b}\right)_{j,k=0}^{n}\\
= & \frac{(a^{2}+bc)^{\frac{n(n+1)}{2}}{\displaystyle \prod_{k=0}^{n}}\left(c|\alpha_{k}|^{2}-2a\Re(\alpha_{k})-b\right){\displaystyle \prod_{j=0}^{k-1}}|\alpha_{k}-\alpha_{j}|^{2}}{{\displaystyle \prod_{k=0}^{n}}{\displaystyle \prod_{j=0}^{k}}|c\alpha_{k}\overline{\alpha_{j}}-a(\alpha_{k}+\overline{\alpha_{j}})-b|^{2}}.\end{align*}
Under the conditions

\begin{align*}
a^{2}+bc\neq0,\, & c\alpha_{j}\overline{\alpha_{k}}-a(\alpha_{j}+\overline{\alpha_{k}})-b\neq0,\,\alpha_{j}\neq\alpha_{k}\end{align*}
for $j,k=0,\dots,n$, the matrix \eqref{eq:generalized muntz matrix}
is invertible, and its inverse $(\gamma_{j,k})_{j,k=0}^{n}$ has element\begin{align*}
\gamma_{j,k} & =\sum_{m=\max(j,k)}^{n}\frac{(c|\alpha_{m}|^{2}-2a\Re(\alpha_{m})-b){\displaystyle \prod_{r=0}^{m-1}}(c\alpha_{r}\overline{\alpha_{j}}-a(\alpha_{r}+\overline{\alpha_{j}})-b)(c\overline{\alpha_{r}}\alpha_{k}-a(\overline{\alpha_{r}}+\alpha_{k})-b)}{(a^{2}+bc)^{m}{\displaystyle \prod_{\begin{array}{c}
p=0\\
p\neq j\end{array}}^{m-1}}(\overline{\alpha_{p}}-\overline{\alpha_{j}}){\displaystyle \prod_{\begin{array}{c}
q=0\\
p\neq k\end{array}}^{m-1}}(\alpha_{q}-\alpha_{k})}.\end{align*}

\end{thm}
When the matrix \eqref{eq:generalized muntz matrix} is positive definite,
its smallest eigenvalue has a lower bound \begin{align*}
\lambda_{s} & \ge\left\{ \sum_{\ell=0}^{n}\left\{ \sum_{j=0}^{\ell}\frac{\sqrt{c|\alpha_{\ell}|^{2}-2a\Re(\alpha_{\ell})-b}{\displaystyle \prod_{j=0}^{\ell-1}}|c\alpha_{j}\overline{\alpha_{k}}-a(\alpha_{j}+\overline{\alpha_{k}})-b|}{(a^{2}+bc)^{\ell/2}{\displaystyle \prod_{\begin{array}{c}
k=0\\
k\neq j\end{array}}^{\ell-1}}|\alpha_{j}-\alpha_{k}|}\right\} ^{2}\right\} ^{-1}.\end{align*}

\subsection{Matrices Associated with $q$-Orthogonal polynomials}

Recall that for $a\in\mathbb{C}$ and $q\in(0,1)$, \cite{Ismail},

\begin{align*}
(a;q)_{\infty} & =\prod_{m=0}^{\infty}(1-aq^{m}),\end{align*}
 \begin{align*}
(a;q)_{m} & =\frac{(a;q)_{\infty}}{(aq^{m};q)_{\infty}},\quad m\in\mathbb{Z},\end{align*}

\begin{align*}
\left[\begin{array}{c}
m\\
j\end{array}\right]_{q} & =\frac{(q;q)_{m}}{(q;q)_{j}(q;q)_{m-j}},\quad0\le j\le m,\end{align*}
 and

\begin{align*}
(a_{1},a_{2},...,a_{n};q)_{m} & =\prod_{k=1}^{n}(a_{k};q)_{m},\quad m\in\mathbb{Z},n\in\mathbb{N}\end{align*}
 for $a_{1},a_{2},...,a_{n}\in\mathbb{C}$. 

For $\nu>-1$, let 

\begin{align*}
0 & <j_{\nu,1}(q)<j_{\nu,2}(q)<\dots<j_{\nu,n}(q)<\dots\end{align*}
be the positive zeros of $z^{-\nu}J_{\nu}^{(2)}(z;q)$, where the
Jackson's $q$-Bessel function $J_{\nu}^{(2)}(z;q)$ is defined as\begin{align*}
J_{\nu}^{(2)}(z;q): & =\frac{(q^{\nu+1};q)_{\infty}}{(q;q)_{\infty}}\sum_{n=0}^{\infty}\frac{(-1)^{n}q^{n(\nu+n)}}{(q,q^{\nu+1};q)_{n}}\left(\frac{z}{2}\right)^{\nu+2n}.\end{align*}
Let us define\begin{align*}
s_{n,\nu} & =\sum_{k=1}^{\infty}\frac{w_{\ell}}{j_{\nu,\ell}^{2n}(q)},\end{align*}
where\begin{align*}
w_{\ell} & =\frac{-4J_{\nu+1}^{(2)}\left(j_{\nu,\ell}(q);q\right)}{\partial_{z}J_{\nu}^{(2)}\left(z;q\right)}\left|_{z=j_{\nu,\ell}(q)}\right.,\end{align*}
 we have the following result:
\begin{thm}
\label{thm:lommel}For $\nu>-1$ and $n\in\mathbb{N}$, the matrix\begin{align}
 & \left(s_{j+k+1,\nu}\right)_{j,k=0}^{n}\label{eq:lommel matrix}\end{align}
has determinant\begin{align*}
\det\left(s_{j+k+1,\nu}\right)_{j,k=0}^{n} & =\frac{2^{-n(n+1)}q^{n(n+1)(4n+6\nu+5)}}{(q^{\nu+1};q^{2})_{n+1}{\displaystyle \prod_{m=1}^{n}(q^{\nu+1};q)_{2m}}},\end{align*}
and its inverse matrix $\left(\gamma_{j,k}\right)_{j,k=0}^{n}$ has
element\begin{align*}
\gamma_{j,k} & =(-4)^{j+k}q^{j(j-\nu)+k(k-\nu)}\sum_{\ell=\max(j,k)}^{n}\left[\begin{array}{c}
\ell+j\\
\ell-j\end{array}\right]_{q}\left[\begin{array}{c}
\ell+k\\
\ell-k\end{array}\right]_{q}\\
\times & \left\{ \frac{(1-q^{2\ell+\nu+1})}{q^{(2j+2k+1)\ell}}\frac{(q^{\nu+1};q)_{\ell+j}}{(q^{\nu+1};q)_{\ell-j}}\frac{(q^{\nu+1};q)_{\ell+k}}{(q^{\nu+1};q)_{\ell-k}}\right\} .\end{align*}
 The smallest eigenvalue of the matrix \eqref{eq:lommel matrix} has
a lower bound
\end{thm}
\begin{align*}
\lambda_{s} & \ge\left\{ \sum_{\ell=0}^{n}\frac{(1-q^{2\ell+\nu+1})h_{2\ell,\nu+1}^{2}(i;q)}{q^{2\ell\nu+\ell(2n+1)}}\right\} ^{-1},\end{align*}
 where \begin{align*}
h_{2n,\nu}(x;q) & =q^{n(n+\nu-1)}\sum_{k=0}^{n}\frac{(q^{\nu},q;q)_{n+k}(-4x^{2})^{k}q^{k(k-2n-\nu+1)}}{(-1)^{n}(q^{\nu},q;q)_{n-k}(q;q)_{2k}}.\end{align*}
Matrices associated with Askey-Wilson orthogonal polynomials:
\begin{thm}
\label{thm:askey}For $n\in\mathbb{N}$, the matrix\begin{align}
 & \left(\frac{(\alpha;q)_{j+k}}{(\alpha\beta;q)_{j+k}}\right)_{j,k=0}^{n}\label{eq:askey-wilson matrix}\end{align}
has determinant\begin{align*}
\det\left(\frac{(\alpha;q)_{j+k}}{(\alpha\beta;q)_{j+k}}\right)_{j,k=0}^{n} & =\frac{\alpha^{n(n+1)/2}q^{n(n^{2}-1)/3}{\displaystyle \prod_{m=1}^{n}}(q,\alpha,\beta;q)_{m}}{{\displaystyle \prod_{m=1}^{n}}(\alpha\beta q^{m-1};q)_{m}(\alpha\beta;q)_{2m}}\end{align*}
Under the conditions

\begin{align*}
\alpha\neq0,\, q\neq0,\, q,\alpha,\beta\neq q^{-k},\, k & =0,\dots n,\end{align*}
the matrix \eqref{eq:askey-wilson matrix} is invertible, and its
inverse matrix $(\gamma_{j,k})_{j,k=0}^{n}$ has element

\begin{align*}
\gamma_{j,k} & =\frac{(-1)^{j+k}q^{\binom{j+1}{2}+\binom{k+1}{2}}}{(\alpha;q)_{j}(\alpha;q)_{k}}\sum_{m=0}^{n}\frac{\left[\begin{array}{c}
m\\
j\end{array}\right]_{q}\left[\begin{array}{c}
m\\
k\end{array}\right]_{q}(\alpha\beta q^{m-1};q)_{j}(\alpha\beta q^{m-1};q)_{k}}{(\alpha q^{j+k})^{m}(q,\beta;q)_{m}(\alpha\beta;q)_{2m}(\alpha\beta q^{m-1};q)_{m}}.\end{align*}
When \textup{$0<q,\alpha,\beta<1$} the matrix \eqref{eq:askey-wilson matrix}
is positive definite, its smallest eigenvalue of the matrix has a
lower bound \begin{align*}
\lambda_{s} & \ge\frac{(\alpha\beta;q)_{\infty}}{(\alpha;q)_{\infty}}\left\{ \sum_{\ell=0}^{n}\frac{p_{\ell}^{2}(-1;\alpha q^{-1},\beta q^{-1}|q)}{h_{\ell}(\alpha q^{-1},\beta q^{-1}|q)}\right\} ^{-1},\end{align*}
 where\begin{align*}
p_{n}(x;\alpha q^{-1},\beta q^{-1}|q) & ={}_{2}\phi_{1}\left(\begin{array}{c}
q^{-n},\alpha\beta q^{n-1}\\
\alpha\end{array};q;qx\right)\end{align*}
and\begin{align*}
h_{n}(\alpha q^{-1},\beta q^{-1}|q) & =\frac{(\alpha\beta;q)_{\infty}}{(\alpha;q)_{\infty}}\frac{(1-\alpha\beta q^{-1})(q,\beta;q)_{n}\alpha^{n}}{(1-\alpha\beta q^{2n-1})(\alpha,a\beta q^{-1};q)_{n}}.\end{align*}
\end{thm}
\begin{rem}
One could get their classical counterparts for the formulas in Theorem
\ref{thm:lommel} and Theorem \ref{thm:askey} by passing the limit
$q\to1^{-}$ with proper normalizations.
\end{rem}

\section{Proofs }

In this section we prove our results under the same restrictions for
the orthogonal polynomials. But the results clearly hold for more
general cases since all the expressions involved except the third
example are rational functions in their parameters.

\subsection{Proof for Theorem \ref{thm:muntz1}}

Given a sequence of distinct complex numbers $\left\{ \alpha_{k}\right\} _{k=0}^{\infty}$
satisfying $\Re(\alpha_{k})>-\frac{1}{2},$ the orthogonal M\" untz-Legendre
polynomials are defined as \cite{Borwein} \begin{align*}
L_{n}(x;\alpha_{0},\dots,\alpha_{n}) & =\sum_{k=0}^{n}c_{n,k}x^{\alpha_{k}},\ c_{n,k}=\frac{{\displaystyle \prod_{j=0}^{n-1}(\alpha_{k}+\overline{\alpha_{j}}+1)}}{{\displaystyle \prod_{\begin{array}{c}
j=0\\
j\neq k\end{array}}^{n-1}}(\alpha_{k}-\alpha_{j})},\end{align*}
for $n\in\mathbb{N}$ and\begin{align*}
L_{0}(x;\alpha_{0},\dots,\alpha_{n}) & =x^{\alpha_{0}}\end{align*}
They satisfy the following orthogonal relation\begin{align*}
\int_{0}^{1}L_{n}(x;\alpha_{0},\dots,\alpha_{n})\overline{L_{m}(x;\alpha_{0},\dots,\alpha_{m})}dx & =\frac{\delta_{m,n}}{1+2\Re(\alpha_{n})}\end{align*}
for $n,m\in\left\{ 0\right\} \cup\mathbb{N}$. We take \begin{align*}
u_{n}(x)= & v_{n}(x)=x^{\alpha_{n}},\end{align*}
and\begin{align*}
I_{j,k} & =\frac{1}{\alpha_{j}+\overline{\alpha_{k}}+1},\ c_{j,k}=\delta_{j,k}.\end{align*}
The orthonormal polynomials are given by

\begin{align*}
p_{n}(x)= & \sqrt{1+2\Re(\alpha_{n})}L_{n}(x;\alpha_{0},\dots,\alpha_{n}).\end{align*}
Then,\begin{align*}
a_{n,k} & =b_{n,k}=\frac{\sqrt{1+2\Re(\alpha_{n})}{\displaystyle \prod_{j=0}^{n-1}(\alpha_{k}+\overline{\alpha_{j}}+1)}}{{\displaystyle \prod_{\begin{array}{c}
j=0\\
j\neq k\end{array}}^{n-1}}(\alpha_{k}-\alpha_{j})}.\end{align*}
 Theorem \ref{thm:muntz1} follows from Theorem \ref{thm:main} and
Theorem \ref{thm:bounds}.

\subsection{Proof for Theorem \ref{thm:muntz2}}

Assume that \begin{align*}
a^{2}+bc & >0,\end{align*}
and \begin{align*}
c\alpha_{j}\overline{\alpha_{k}}-a(\alpha_{j}+\overline{\alpha_{k}})-b & >0,\end{align*}
the following orthogonal M\" untz polynomials exist, \cite{Marinkovic}
\begin{align*}
q_{n}(x) & =\sum_{k=0}^{n}A_{n,k}x^{\alpha_{k}},\ A_{n,k}=\frac{{\displaystyle \prod_{j=0}^{n-1}}\left(\alpha_{k}-\frac{a\overline{\alpha_{j}}+b}{c\overline{\alpha_{j}}-a}\right)}{{\displaystyle \prod_{\begin{array}{c}
j=0\\
j\neq k\end{array}}^{n-1}}(\alpha_{k}-\alpha_{j})}.\end{align*}
They satisfy the orthogonal relation\begin{align*}
(q_{n}(x),q_{m}(x))_{*} & =h_{n}\delta_{m,n},\end{align*}
where\begin{align*}
h_{n} & =\frac{(a^{2}+bc)^{n}}{\left(c|\alpha_{n}|^{2}-2a\Re(\alpha_{n})-b\right){\displaystyle \prod_{j=0}^{n-1}}|c\alpha_{j}-a|^{2}}.\end{align*}
Take \begin{align*}
u_{n}(x) & =v_{n}(x)=x^{\alpha_{n}},\, c_{j,k}=\delta_{j,k}.\end{align*}
From \cite{Marinkovic} we have \begin{align*}
I_{j,k} & =\frac{1}{c\alpha_{j}\overline{\alpha_{k}}-a(\alpha_{j}+\overline{\alpha_{k}})-b}.\end{align*}
The orthonormal polynomials are\begin{align*}
p_{n}(x) & =\sum_{k=0}^{n}a_{n,k}x^{\alpha_{k}},\end{align*}
where\begin{align*}
a_{n,k} & =b_{n,k}=\frac{A_{n,k}}{\sqrt{h_{n}}}.\end{align*}
Then, Theorem \ref{thm:muntz2} follows from Theorem \ref{thm:main}
and Theorem \ref{thm:bounds}.

\subsection{Proof for Theorem \ref{thm:lommel}}

The even $q$-Lommel polynomials satisfy the orthogonal relation \cite{Ismail}\begin{align*}
\sum_{\ell=1}^{\infty}h_{2n,\nu+1}\left(\frac{1}{j_{\nu,\ell}(q)};q\right)h_{2m,\nu+1}\left(\frac{1}{j_{\nu,\ell}(q)};q\right)\frac{w_{\ell}}{j_{\nu,\ell}^{2}(q)} & =\frac{q^{2n\nu+n(2n+1)}\delta_{m,n}}{1-q^{2n+\nu+1}}.\end{align*}
Let \begin{align*}
u_{j} & =v_{j}=x^{2j},\end{align*}
then\begin{align*}
I_{j,k} & =s_{j+k+1,\nu},\quad c_{j,k}=\delta_{j,k}.\end{align*}
The orthonormal polynomials are given by\begin{align*}
p_{n}(x) & =\frac{\sqrt{1-q^{2n+\nu+1}}}{q^{n\nu+n(n+1/2)}}h_{2n,\nu+1}(x;q).\end{align*}
hence\begin{align*}
a_{n,k} & =b_{n,k}=\frac{\sqrt{1-q^{2n+\nu+1}}(q^{\nu+1},q;q)_{n+k}4^{k}q^{k^{2}-k(2n+\nu)}}{(-1)^{n-k}(q^{\nu+1},q;q)_{n-k}(q;q)_{2k}q^{n/2}}.\end{align*}
Then, Theorem \ref{thm:lommel} follows from Theorem \ref{thm:main}
and Theorem \ref{thm:bounds}.

\subsection{Proof for Theorem \ref{thm:askey} }

We present two proofs to this example. Our first proof uses the Askey-Wilson
polynomials, while the second uses the little $q$-Jacobi polynomials.
For some mysterious reasons, these polynomials yield essentially the
same matrices.

For each nonnegative integer $n$, the Askey-Wilson polynomial has
the following series representation\begin{align*}
a_{n}(x;{\normalcolor t}|q) & =t_{1}^{-n}(t_{1}t_{2},t_{1}t_{3},t_{1}t_{4};q)_{n}\\
\times & {}_{4}\phi_{3}\left(\begin{array}{c}
q^{-n},t_{1}t_{2}t_{3}t_{4}q^{n-1},t_{1}e^{i\theta},t_{1}e^{-i\theta}\\
t_{1}t_{2},t_{1}t_{3},t_{1}t_{4}\end{array}\vert q,q\right),\quad x=\cos\theta,\end{align*}
where the basic hypergeometric function ${}_{r}\phi_{s}$ with complex
parameters $a_{1},...,a_{r};b_{1},...,b_{s}$ is formally defined
as, \cite{Ismail}\begin{align*}
{}_{r}\phi_{s}\left(\begin{array}{c}
a_{1},...,a_{r}\\
b_{1},...,b_{s}\end{array};q,z\right) & =\sum_{n=0}^{\infty}\frac{(a_{1},...,a_{r};q)_{n}z^{n}}{(q,b_{1},...,b_{s};q)_{n}}\left((-1)^{n}q^{(n-1)n/2}\right)^{s+1-r}.\end{align*}
It is well known that $a_{n}(x;{\normalcolor t}|q)$ is symmetric
in the real parameters $t_{1},t_{2},t_{3},t_{4}$. Under the condition
$\max\left\{ |t_{1}|,|t_{2}|,|t_{3}|,|t_{4}|\right\} <1$, the Askey-Wilson
polynomials satisfy the following orthogonal relation\begin{align*}
\int_{-1}^{1}a_{m}(x;{\normalcolor t}|q)a_{n}(x;{\normalcolor t}|q)w(x;{\normalcolor t}|q)dx & =h_{n}\delta_{m,n},\end{align*}
where\begin{align*}
h_{n} & =\frac{2\pi(t_{1}t_{2}t_{3}t_{4}q^{2n};q)_{\infty}(t_{1}t_{2}t_{3}t_{4}q^{n-1};q)_{n}}{(q^{n+1};q)_{\infty}{\displaystyle \prod_{1\le j<k\le4}(t_{j}t_{k}q^{n};q)_{\infty}}}\end{align*}
and\begin{align*}
w(x;{\normalcolor t}|q) & =\frac{(e^{2i\theta},e^{-2i\theta};q)_{\infty}}{{\displaystyle \prod_{j=1}^{4}(t_{j}e^{i\theta},t_{j}e^{-i\theta};q)_{\infty}}}\frac{1}{\sqrt{1-x^{2}}},\quad x=\cos\theta.\end{align*}
Let\begin{align*}
u_{j}(x) & =(t_{1}e^{i\theta},t_{1}e^{-i\theta};q)_{j},\quad v_{j}(x)=(t_{2}e^{i\theta},t_{2}e^{-i\theta};q)_{j},\end{align*}
it is known that \cite{Ismail}\begin{align*}
\frac{(be^{i\theta},be^{-i\theta};q)_{n}}{(ab,b/a;q)_{n}} & =\sum_{k=0}^{n}\frac{(q^{-n},ae^{i\theta},ae^{-i\theta})_{k}q^{k}}{(q,ab,q^{1-n}a/b;q)_{k}},\end{align*}
for $a\cdot b\neq0$, hence,\begin{align*}
c_{n,k} & =\left[\begin{array}{c}
n\\
k\end{array}\right]_{q}\left(t_{1}t_{2}q^{k};q\right)_{n-k}\left(\frac{t_{2}}{t_{1}};q\right)_{n-k}\left(\frac{t_{2}}{t_{1}}\right)^{k}.\end{align*}
and

\begin{align*}
I_{j,k} & =\int_{-1}^{1}w(x;t_{1}q^{j},t_{2}q^{k},t_{3},t_{4};q)dx,\end{align*}
or

\begin{align*}
I_{j,k} & =\frac{2\pi(t_{1}t_{2}t_{3}t_{4}q^{j+k};q)_{\infty}}{(q,t_{1}t_{2}q^{j+k},t_{1}t_{3}q^{j},t_{1}t_{4}q^{j},t_{2}t_{3}q^{k},t_{2}t_{4}q^{k},t_{3}t_{4};q)_{\infty}}.\end{align*}
The orthonormal polynomials are

\begin{align*}
p_{j}(x) & =\frac{a_{j}(x;{\normalcolor t}|q)}{\sqrt{h_{j}}},\end{align*}
with\begin{align*}
a_{n,k} & =\frac{(t_{1}t_{2},t_{1}t_{3},t_{1}t_{4};q)_{n}(q^{-n},t_{1}t_{2}t_{3}t_{4}q^{n-1};q)_{k}q^{k}}{(t_{1})^{n}\sqrt{h_{n}}(q,t_{1}t_{2},t_{1}t_{3},t_{1}t_{4};q)_{k}},\end{align*}
 and\begin{align*}
b_{n,k} & =\frac{(t_{2}t_{1},t_{2}t_{3},t_{2}t_{4};q)_{n}(q^{-n},t_{1}t_{2}t_{3}t_{4}q^{n-1};q)_{k}q^{k}}{(t_{2})^{n}\sqrt{h_{n}}(q,t_{2}t_{1},t_{2}t_{3},t_{2}t_{4};q)_{k}}.\end{align*}
 Then\begin{align*}
\det H_{n} & =\prod_{m=0}^{n}\frac{2\pi\left(t_{1}t_{2}q^{m-1}\right)^{m}(t_{1}t_{2}t_{3}t_{4}q^{2m};q)_{\infty}}{(t_{1}t_{2}t_{3}t_{4}q^{m-1};q)_{m}(q^{m+1},q)_{\infty}{\displaystyle {\displaystyle {\displaystyle {\displaystyle \prod_{1\le j<k\le4}}}}}(t_{j}t_{k}q^{m};q)_{\infty}},\end{align*}
 which could be simplified to \begin{align}
 & \det\left(\frac{(t_{1}t_{2};q)_{j+k}}{(t_{1}t_{2}t_{3}t_{4};q)_{j+k}}\right)_{j,k=0}^{n}\label{eq:askey1}\\
= & \prod_{m=0}^{n}\frac{\left(t_{1}t_{2}q^{m-1}\right)^{m}(q,t_{1}t_{2},t_{3}t_{4};q)_{m}}{(t_{1}t_{2}t_{3}t_{4};q)_{2m}(t_{1}t_{2}t_{3}t_{4}q^{m-1};q)_{m}}.\nonumber \end{align}
For any $n=0,1,...,$ let $e,c_{0},c_{1},...,c_{n},d_{0},d_{1},...,d_{n}$
be non-zero numbers and, \begin{align*}
X & =\left(x_{jk}\right)_{j,k=0}^{n},\ Y=\left(y_{jk}\right)_{j,k=0}^{n},\ I=\left(\delta_{j,k}\right)_{j,k=0}^{n}.\end{align*}
If \begin{align*}
XY & =I,\end{align*}
then,

\begin{align*}
\widetilde{X}\widetilde{Y} & =I,\end{align*}
where \begin{align*}
\widetilde{X} & =\left(ex_{jk}c_{j}d_{k}\right)_{j,k=0}^{n},\end{align*}
 and \begin{align*}
\widetilde{Y} & =\left(\frac{y_{jk}}{ed_{j}c_{k}}\right)_{j,k=0}^{n}.\end{align*}
Using the above trick, we simplify the inverse pairs down to\begin{align*}
\left(\frac{(t_{1}t_{2};q)_{j+k}}{(t_{1}t_{2}t_{3}t_{4};q)_{j+k}}\right)_{j,k=0}^{n}\end{align*}
and its inverse matrix $(\omega_{j,k})_{j,k=0}^{n}$ with element
\begin{align}
\omega_{j,k} & =\frac{(-1)^{j+k}q^{\binom{j+1}{2}+\binom{k+1}{2}}}{(t_{1}t_{2};q)_{j}(t_{1}t_{2};q)_{k}}\sum_{m=0}^{n}\frac{\left[\begin{array}{c}
m\\
j\end{array}\right]_{q}\left[\begin{array}{c}
m\\
k\end{array}\right]_{q}}{(q,t_{3}t_{4};q)_{m}}\label{eq:askey2}\\
\times & \left\{ \frac{(t_{1}t_{2}t_{3}t_{4}q^{m-1};q)_{j}(t_{1}t_{2}t_{3}t_{4}q^{m-1};q)_{k}}{\left(t_{1}t_{2}q^{j+k}\right)^{m}(t_{1}t_{2}t_{3}t_{4};q)_{2m}(t_{1}t_{2}t_{3}t_{4}q^{m-1};q)_{m}}\right\} \nonumber \end{align}
Except the lower bound, the assertions of Theorem \ref{thm:askey}
follows from \eqref{eq:askey1} and \eqref{eq:askey2} by the change
of variables, \begin{align*}
\alpha & =t_{1}t_{2},\quad\beta=t_{3}t_{4}.\end{align*}

Assume that \begin{align*}
p_{-1}(x;a,b|q) & =0,\quad p_{0}(x;a,b|q)=1,\end{align*}
the little $q$-Jacobi polynomials $\left\{ p_{n}(x;a,b|q)\right\} _{n=0}^{\infty}$
have the orthogonal relation \cite{Ismail}

\begin{align*}
\sum_{k=0}^{\infty}\frac{(bq;q)_{k}(aq)^{k}}{(q;q)_{k}}p_{m}(q^{k};a,b|q)p_{n}(q^{k};a,b|q) & =h_{n}(a,b|q)\delta_{mn}\end{align*}
for $m,n\ge0$ The moments are given by the formula\begin{align*}
\mu_{n} & =\sum_{m=0}^{\infty}\frac{(bq;q)_{m}(aq)^{m}q^{nm}}{(q;q)_{m}},\end{align*}
 or\[
\mu_{n}=\frac{(abq^{n+2};q)_{\infty}}{(aq^{n+1};q)_{\infty}},\]
 by using the $q$-binomial theorem \cite{Ismail}. The orthonormal
polynomial\begin{align*}
p_{n}(x) & =\frac{(-1)^{n}p_{n}(x;aq,bq|q)}{\sqrt{h_{n}(a,b|q)}}.\end{align*}
Hence, \begin{align*}
a_{n,k} & =b_{n,k}=\frac{(-1)^{n}(q^{-n},abq^{n+1};q)_{k}q^{k}}{\sqrt{h_{n}(a,b|q)}(q,aq;q)_{k}}.\end{align*}
Theorem \ref{thm:askey} follows from Theorem \ref{thm:main} and
Theorem \ref{thm:bounds} with a change of variables \begin{align*}
\alpha & =aq,\quad\beta=bq.\end{align*}

\begin{acknowledgement*}
This work is partially supported by Chinese National Natural Science
Foundation grant No.10761002, Guangxi Natural Science Foundation grant
No.0728090.\end{acknowledgement*}

\end{document}